\documentclass[10pt]{article}

\usepackage{a4wide}
\usepackage{amssymb}
\usepackage{amsfonts}
\usepackage{enumitem}
\usepackage{amsmath}
\usepackage[all]{xy}
\input xy

\date{}

 \newtheorem{proposition}{Proposition}[section]
\newtheorem{theorem}[proposition]{Theorem}

\newtheorem{definition}[proposition]{Definition}
\newtheorem{corollary}[proposition]{Corollary}

\def\der{\partial }

\def\nFM0{{\nu }_{F,M_0}}
\def\nFN0{{\nu }_{F,N_0}}
\def\nGN0{{\nu }_{G,N_0}}

\def\N0{ {\bf N}_0 }

\def\ra{\rightarrow}

\def\Xpm{X^{\pm }}

\def\l1{{\lambda}_1}

\def\a{\alpha}
\def\a0{ {\alpha }_0}
\def\a1{ {\alpha }_1}

\def\l{\lambda}

\def\nFGM0{{\nu }_{F,G,M_0}}

\def\nFN0{{\nu}_{F,N_0}}

\def\sm{{\sigma}^m}

\def\sm1{{\sigma}^{-1}}

\def\smtp1{{\sigma}^{-t+1}}

\def\S1{S^{-1}}

\def\Xpm1{X^{\pm 1}_1}

\def\sPM1{{\sigma }^{\pm 1}}
\def\sMP1{{\sigma }^{\mp 1 }}

\def\di{{\rm d.ind}}

\def\L{\Lambda}

\def\Ytm1{Y^{t-1}}
\def\Yim1{Y^{i-1}}

\def\CM{{\cal M}}

\def\CF{{\cal F}}

\def\dim{{\rm dim }}

\def\SL2Z{ {\rm SL}_2({\bf Z}) }

\def\th{ \theta }

\def\Gp1{ G^{1 , 1 } }
\def\P11{ P^{-1 , 1 } }
\def\Pp1{ P^{1 , 1 } }

\def\th{\theta}

\def\nCLsr{{}^\nu\kern-2pt {\cal L}^{\sigma , \rho  }}
\def\nP{{}^\nu \kern-2pt P}
\def\nL{{}^\nu\kern-2pt L}
\def\nLL{{}^\nu\kern-2pt \Lambda}
\def\nPsr{{}^\nu\kern-2pt P^{\sigma , \rho  }}
\def\nLsr{{}^\nu\kern-2pt L^{\sigma , \rho  }}
\def\nuCL{{}^\nu\kern-2pt  {\cal L}}
\def\nCLsr{{}^\nu\kern-2pt {\cal L}^{\sigma , \rho  }}
\def\nCL1m{{}^\nu\kern-2pt {\cal L}^{-1 , 1  }}
\def\x1nu{x^\frac{1}{\nu}}
\def\xm1nu{x^{-\frac{1}{\nu}}}

\def\ra{\rightarrow }

\def\CB{{\cal B}}

\def\CI{{\cal I}}

\def\CT{{\cal T}}

\def\nAM0{{\nu }_{{\cal A},M_0}}
\def\nAN0{{\nu }_{{\cal A},N_0}}

\def\End{ {\rm End }}

\def\ga{\mathfrak{a}}

\def\gm{\mathfrak{m}}
\def\gp{\mathfrak{p}}

\def\SL{{\rm SL}}

\def\di!{\frac{\der^i}{i!}}
\def\dik!{\frac{\der^k_i}{k!}}

\def\N{\mathbb{N}}

\def\0{\overline{0}}
\def\1{\overline{1}}

\def\Ln1{\L_{n,\overline{1}}}

\def\a1{a_{\overline{1}}}

\def\S{\Sigma}

\def\vn1{\overrightarrow{n-1}}

\def\Sub{{\rm Sub}}

\def\mJ{\mathbb{J}}
\def\mI{\mathbb{I}}

\def\ann{{\rm ann}}

\def\K1{{\rm K}_1}

\def\hmI1{\widehat{\mI_1}}
\def\tmI1{\widetilde{\mI_1}}
\def\tmJ1{\widetilde{\mJ_1}}
\def\hB1{\widehat{B_1}}
\def\hCB1{\widehat{\CB_1}}

\def\ga{\mathfrak{a}}

\def\Cyc{{\rm Cyc}}

\def\ann{{\rm ann}}

\def\card{{\rm card}}

\def\Sub{{\rm Sub}}

\begin{document}

\author{ V. V. \ Bavula 
}

\title{Classification of multiplication modules over multiplication  rings with finitely many minimal primes}


\maketitle

\begin{abstract}
A classification of multiplication modules over multiplication  rings with finitely many minimal primes is obtained.  A characterisation of multiplication rings with finitely many minimal primes is given via faithful, Noetherian, distributive modules. It is proven that  for a multiplication  ring with finitely many minimal primes  every faithful, Noetherian, distributive module is a faithful multiplication module, and vice versa.

$\noindent $






\end{abstract}



\section{Introduction}\label{INTR}

In this paper, all rings are commutative with 1 and all modules are unital. A ring $R$ is called a {\bf multiplication ring} if $I$ and $J$ are ideals of $R$ such that $J\subseteq I$ then $J=I^{\prime}I$ for some ideal $I^{\prime}$ of $R$. An $R$-module $M$ is called a {\bf multiplication module} if each submodule of $M$ is equal to $IM$ for some ideal $I$ of the ring $R$. The concept of multiplication ring was introduced by Krull in \cite{Krull-1948}. In \cite{Mott-1969}, Mott proved that a multiplication ring has finitely many minimal prime ideals iff it is a Noetherian ring.\\

The next theorem is a description of multiplication rings with finitely many minimal primes.

\begin{theorem}\label{C4May18}
(\cite[Theorem 1.1]{Als-Bav-N1}) Let $R$ be a ring with finitely many minimal prime ideals. Then the ring $R$ is a multiplication ring iff $R\cong \displaystyle \prod_{i=1}^{n} D_{i}$ is a finite direct product of rings where $D_{i}$ is either a Dedekind domain or an Artinian, local principal ideal ring.
\end{theorem}

{\bf Classification  of multiplication modules over multiplication  rings with finitely many minimal primes.} 
Using Theorem \ref{C4May18}, a  criterion for a direct sum of modules to be a multiplication module (Theorem \ref{1.20})  and some other results, a classification of multiplication modules over a multiplication  ring with finitely many minimal primes is given, Theorem  \ref{2Dec18}.
 
\begin{theorem}\label{2Dec18}
 Let $R$ be a multiplication  ring with finitely many minimal primes, i.e., 
  $R\cong \displaystyle \prod_{i=1}^{n} D_{i}$ is a finite direct product of rings where $D_{i}$ is either a Dedekind domain or an Artinian, local principal ideal ring and $1=e_1+\cdots + e_n$ be the corresponding sum of central orthogonal idempotents of the ring $R$. Let $M$ be an $R$-modules and $M=\oplus_{i=1}^n M_i$ where $M_i:= e_iM$. 
 Then the $R$-module $M$ is a multiplication $R$-module iff each $D_i$-module $M_i$ is either isomorphic to  $D_i$ or to $D_i/I_i$ where $I_i$ is a nonzero ideal of $D_i$ or to  a nonzero ideal of the ring $D_i$ in case when the ring $D_i$ is a Dedekind domain. 
\end{theorem}

{\bf Classification  of faithful multiplication modules over a multiplication  ring with finitely many minimal primes.}

\begin{theorem}\label{A9Dec18}
 Let $R$ be a multiplication  ring with finitely many minimal primes. We keep the notation of Theorem \ref{2Dec18} ($R\cong \displaystyle \prod_{i=1}^{n} D_{i}$).  Then an $R$-module $M=\oplus_{i=1}^n M_i$ (where $M_i= e_iM$) is a faithful multiplication $R$-module iff for each $i=1, \ldots , n$, either ${}_RM_i\simeq D_i$ or  ${}_RM_i\simeq I_i$ where $I_i$ is a nonzero ideal of 
 the ring $D_i$  in case when $D_i$ is a Dedekind domain. 
\end{theorem}

{\it Proof.} The theorem follows at once from Theorem \ref{2Dec18}. $\Box $\\

{\bf Characterisation of multiplication rings with finitely many minimal primes via faithful, Noetherian, distributive modules.} Let $R$ be a ring and $M$ be an $R$-module. A submodule $N$ of $M$ is called a {\em distributive submodule} if one of the following equivalent conditions holds:
\begin{eqnarray*}
(M_1+M_2)\cap N  &=& M_1\cap N +M_2\cap N, \\
M_1\cap M_2+ N &=&(M_1+N)\cap (M_2+N).
\end{eqnarray*}
The $R$-module $M$ is called a {\bf distributive module} if all submodules of $M$ are distributive submodules.

\begin{theorem}\label{9Dec18}
A commutative ring $R$ is a multiplication ring with finitely many minimal primes iff there is a faithful, Noetherian, distributive $R$-module.
\end{theorem}

{\bf Classification  of faithful, Noetherian, distributive modules   over a multiplication  ring with finitely many minimal primes.}

\begin{theorem}\label{B9Dec18}
 Let $R$ be a multiplication  ring with finitely many minimal primes. Then every faithful, Noetherian, distributive $R$-module is a faithful multiplication $R$-module, and vice versa.
\end{theorem}


\section{Proofs}\label{PRFS}

In this section we prove the results from the Introduction.

\begin{definition}\label{a1.20}
We say that the {\bf intersection condition} holds for a direct sum $M=\bigoplus_{\lambda \in \Lambda} M_{\lambda}$ of nonzero $R$-modules $M_{\lambda}$ if for all submodules $N$ of $M$, $N= \bigoplus_{\lambda\in \Lambda} (N\bigcap M_\lambda )$.
\end{definition}

\begin{definition}\label{a1.28}
Let $M=\bigoplus_{\lambda\in \Lambda} M_{\lambda}$ be a direct sum of nonzero $R$-modules with $\card(\Lambda )\geqslant 2$, $\ga_{\lambda}=\ann_{R}(M_{\lambda})$ and $\ga^{\prime}_{\lambda}=\cap_{\mu\neq \lambda}\ga_{\mu}$. We say that the {\bf orthogonality condition} holds for the direct sum $M=\bigoplus_{\lambda\in \Lambda} M_{\lambda}$ if $\ga^{\prime}_{\lambda}M_{\mu}= \delta_{\lambda\mu}M_{\mu}$ for all $\lambda,\mu\in \Lambda$. Clearly, $\ga^{\prime}_{\lambda}\neq 0$ for all $\lambda \in \Lambda$ (since all $M_{\lambda} \neq 0$). In particular, $\ga_{\lambda} \neq 0$ for all $\lambda \in \Lambda$.
\end{definition}

\begin{definition}\label{C1.20}
Let $M=\bigoplus_{\lambda \in \Lambda} M_{\lambda}$ be a direct sum of nonzero $R$-modules with $\card(\Lambda) \geq 2$. We say that the {\bf strong orthogonality condition} holds for $M$ if for each set of $R$-modules $\lbrace N_{\lambda} \rbrace_{\lambda\in \Lambda}$ such that $N_{\lambda}\subseteq M_{\lambda}$, there is a set of of ideals $\lbrace I_{\lambda} \rbrace_{\lambda\in \Lambda}$ of $R$ such that $I_{\lambda} M_{\mu}=\delta_{\lambda \mu}N_{\lambda}$ for all $\lambda,$ $\mu\in \Lambda$ where $\delta_{\lambda \mu}$ is the Kronecker delta. The set of ideals $\lbrace I_{\lambda} \rbrace_{\lambda\in \Lambda}$ is called an {\bf orthogonalizer} of $\lbrace N_{\lambda} \rbrace_{\lambda\in \Lambda}$. 
\end{definition}

Theorem \ref{1.20} is one of the  criteria for a direct sum of modules to be a multiplication module that are obtained in \cite{Als-Bav-N1}. It is  given via the intersection and strong orthogonality conditions.

\begin{theorem}\label{1.20}

\cite{Als-Bav-N2} Let $M=\bigoplus_{\lambda \in \Lambda} M_{\lambda}$ be a direct sum of nonzero $R$-modules with $\card(\Lambda) \geq 2$. Then $M$ is a multiplication module iff the intersection and strong orthogonality conditions hold for the direct sum $M=\bigoplus_{\lambda \in \Lambda} M_{\lambda}$.
\end{theorem}

 An $R$-module is called a {\em cyclic} if it is 1-generated. For an $R$-module $M$, let $\Cyc_R(M)$ be the set of its {\em cyclic} submodules. For an $R$-module $M$, we denote by $\ann_{R}(M)$ its annihilator. An $R$-module $M$ is called {\em faithful} if $\ann_{R}(M)=0$. For a submodule $N$ of $M$, the set $[N:M]:= \ann_R(M/N)= \{ r\in R\, | \, rM \subseteq N\}$ is an ideal of the ring $R$ that contains the {\em annihilator} $\ann_R(M)=[0:M]$ of the module $M$. The set $\th(M):=\sum_{C\in \Cyc_R(M)}[C:M]$ is an ideal of $R$. Clearly, $\ann_R(M)\subseteq \th(M)$. If $M$ is an ideal of $R$ then $M\subseteq \theta(M)$.\\

{\bf Proof of Theorem \ref{2Dec18}.} $(\Leftarrow )$ All the $D_i$-modules $M_i$ of the theorem are multiplication  $D_i$-modules. Hence, the direct sum $\oplus_{i=1}^nM_i$ is a multiplication module over the direct product rings $R=\prod_{i=1}^nD_i$. \\

$(\Rightarrow )$   Suppose that the $R$-module $M=\oplus_{i=1}^n M_i$ is a multiplication $R$-module where $M_i=e_iM$ for $i=1, \ldots , n$.\\ 
 
 (i) {\em The $D_i$-module $M_i$ is a multiplication $D_i$-module:} The statement is obvious since $R=\prod_{i=1}^nD_i$. \\
 
 (ii) {\it The $D_i$-module $M_i$ is a finitely generated  $D_i$-module:} Since $M_i$ is a multiplication $D_i$-module, $$M_i=\sum_{C\in \Cyc_{D_i}(M_i)}C=\sum_{C\in \Cyc_{D_i}(M_i)}[C:M_i]M_i=(\sum_{C\in \Cyc_{D_i}(M_i)}[C:M_i])M_i=\th (M_i)M_i.$$ The ideal $\th (M_i)=\sum_{C\in \Cyc_{D_i}(M_i)}[C:M_i]$ of the Noetherian ring $D_i$ is a finitely generated $D_i$-module, i.e., $\th (M_i) = \sum_{i=1}^{n_i}D_i\th_i$  for some elements $\th_i\in \th (M_i)$. Then $$M_i=\th (M_i) M_i=\sum_{i=1}^{n_i}D_i\th_iM_i\subseteq \sum_{i=1}^{n_i}C_i\subseteq M_i,$$ and so the $D_i$-module $M_i=\sum_{i=1}^{n_i} C_i$ is finitely generated. \\

(iii) {\em Suppose that the ring $D_i$ is a Dedekind domain. Then the $D_i$-module $M_i$ is isomorphic either to $D_i$ or to $D_i/I_i$ or to $J_i$ where $I_i$ and $J_i$ are ideals of the ring $D_i$:} It is well-known that a nonzero finitely generated module $\CM$ over a Dedekind domain $D$ is a direct sum $\CM = \CF \oplus \CT$ of a torsionfree $D$-module $\CF$ and a torsion $D$-module $\CT$; $\CF = I\oplus D^m$ for some ideal $I$ of $D$ and $m\geq 0$; and $\CT =\oplus_{i=1}^{t_i}D/\gp_i^{m_i}$ where $\gp_i$ are maximal ideals of the ring $D$ and  $m_i\in \N$. Suppose that the $D$-module $\CM$ is a multiplication $D$-module. By Theorem \ref{1.20}, the direct sum of $D$-modules $$\CM = I\oplus D^m \oplus \bigoplus_{i=1}^{t_i} D/\gp_i^{m_i}$$ must satisfy the strong orthogonality conditions. Hence, either $\CM = I$ of $\CM = D$ or $\CM = \oplus_{i=1}^{t_i} D/\gp_i^{m_i}$ where $\gp_1, \ldots , \gp_{t_i}$ are {\em distinct} maximal ideals of the ring $D$, and so $\CM =\oplus_{i=1}^{t_i} D/\gp_i^{m_i}\simeq D/\prod_{i=1}^{t_i}\gp_i^{m_i}$. \\

(iv) {\em Suppose that $D_i$ is an Artinian, local, principal ideal ring. Then the $D_i$-module $M_i$ is isomorphic  either to $D_i$ or to $D_i/I_i$ where $I_i$ is a nonzero ideal of $D_i$:} Let $D=D_i$ and $\gm  $ be the maximal ideal of the local ring $D_i$ and $\gm^{\nu } \neq 0$ and $\gm^{\nu +1} =0$ for some natural number $\nu$. Then $$\{ D, \gm , \gm^2, \ldots , \gm^\nu , \gm^{\nu +1} =0\}$$  is the set of all the ideals of the ring $D$. The $D$-module $M_i$ is a nonzero finitely generated multiplication $D$-module. Hence, $\{ M_i, \gm M_i, \gm^2 M_i, \ldots , \gm^\mu  M_i, \gm^{\mu +1} M_i =0\}$  is the set of all $D$-submodules of $M_i$ for some natural number $\mu$ such that $\mu \leq \nu$. In particular, the $D$-module $M_i$ is a uniserial $D$-module since
$$M_i\supset  \gm M_i \supset   \gm^2 M_i\supset   \cdots \supset   \gm^\mu  M_i\supset    \gm^{\mu +1} M_i =0.$$
Therefore, $$\dim_{k_\gm}(M_i/\gm M_i)=1$$ where $k_\gm :=D/\gm $, and so $M_i= Dm_i+\gm M_i$ for some element $m_i\in M_i\backslash \gm M_i$. By the Nakayama Lemma, $M_i=D\gm_i$, and the statement (iv) follows. $\Box$\\

\begin{corollary}\label{a2Dec18}
Let $R$ be an Artinian  multiplication  ring. Then every multiplication $R$-module if an epimorphic image of the $R$-module $R$.
\end{corollary}

{\it Proof}. The corollary follows at once from Theorem \ref{2Dec18}. $\Box $

\begin{corollary}\label{b2Dec18}
Let $R$ be a multiplication  ring with finitely many minimal primes and $M$ be a multiplication $R$-module. Then 
\begin{enumerate}
\item The endomorphism ring $\End_R(M)$ is also a multiplication ring.
\item $\End_R(M)\simeq R/\ann_R(M)$.
\item The $\End_R(M)$-module $M$ is a faithful multiplication $\End_R(M)$-module.
\end{enumerate}
\end{corollary}

{\it Proof}. The corollary follows at once from Theorem \ref{2Dec18}. $\Box $ \\

In the proof of Theorem \ref{9Dec18} we will use the following results.

\begin{theorem}\label{RES-T}
Let $R$ be a commutative ring. 
\begin{enumerate}
\item (\cite[Corollary, p.177]{Barnard-MM-1981}) Let $M$ be a Noetherian distributive $R$-module. Then every submodule of $M$ which is locally nonzero at every maximal ideal of $R$, is of the form $IM$ where $I$ is a unique product of maximal ideals of $R$. 
\item (\cite[Lemma 2.(ii)]{Barnard-MM-1981}) A finitely generated $R$-module $M$ is a multiplication module iff the $R_\gp$-module $M_\gp$ is a multiplication module for all prime/maximal ideals $\gp$ of $R$. 
\item (\cite[Theorem 1.3.(ii)]{El-Bast-Smith-MM-1981}) {\sc (Cancellation Law)} If $M$ is a finitely generated, faithful multiplication $R$-module then for any two ideals $A$ and $B$ of $R$, $MA\subseteq BM$ iff $A\subseteq B$. 
\end{enumerate}
\end{theorem}

{\bf Proof of Theorem \ref{9Dec18}.} $(\Rightarrow )$ By Theorem \ref{2Dec18}, the $R$-module $R$ is a faithful, Noetheriaan, distributive $R$-module. 

$(\Leftarrow )$ Let $M$ be  faithful, Noetheriaan, distributive $R$-module. 

(i) {\it The ring $R$ is a Noetherian ring}: The $R$-module $M$ is Noetherian, hence finitely generated, $M=\sum_{i=1}^n Rm_i$ for some elements $m_1, \ldots , m_n\in M$. The $R$-module $M$ is a faithful module. Hence, the map $R\ra \oplus_{i=1}^n Rm_i$, $r\mapsto (rm_1, \ldots , rm_n)$ is an $R$-monomorphism. The direct sum is a Noetherian $R$-module (as a finite direct sum of Noetherian modules),  and the statement (i) follows. 

(ii) {\it The ring $R$ has only finitely many minimal primes}: The statement (ii) follows from the statement (i). 

(ii) {\it For all maximal ideals $\gm$ of the ring $R$, the $R_\gm$-module $M_\gm$ is faithful, Noetherian and distributive}: The $R$-module $M$ is finitely generated. Hence, $\ann_{R_\gm}(M_\gm) = \ann_R(M)_\gm=0$ since $\ann_R(M)=0$. Clearly, the $R_\gm$-module $M_\gm$ is  Noetherian and distributive (since the $R$-module $M$ is so and localizations respect finite intersections). 

(iv) {\it The $R_\gm$-module $M_\gm$ is a multiplication $R_\gm$-module}: The $R$-module $M$ is finitely generated. By the statement (iv) and Theorem \ref{RES-T}.(2), the $R$-module $M$ is a multiplication $R$-module.

Let ($\CI (R), \subseteq )$ be the lattice of ideals of the ring $R$ and $(\Sub_R(M), \subseteq )$ be the lattice of $R$-submodules of the $R$-module $M$. 

(vi) {\it The map $\CI (R)\ra \Sub_R(M)$, $ I\mapsto IM$ is an isomorphism of latices}: The $R$-module $M$ is a finitely generated, faithful multiplication module (the statement (v)), and the statement (vi) follows from Theorem \ref{RES-T}.(3).

(vii) {\it The ring $R$ is a multiplication ring}:  The statement (vii) follows from the statements (v) and (vi).

Now, the theorem follows from the statement (ii) and (vii). $\Box$  \\

{\bf Proof of Theorem \ref{B9Dec18}.} $(\Rightarrow )$ See the statement (vi) in the proof of Theorem \ref{9Dec18}.

$(\Leftarrow )$ This implication follows at once from Theorem \ref{A9Dec18}. $\Box$

{\bf Licence.} For the purpose of open access, the author has applied a Creative Commons Attribution (CC BY) licence to any Author Accepted Manuscript version arising from this submission.

{\bf Disclosure statement.} No potential conflict of interest was reported by the author.

{\bf Data availability statement.} Data sharing not applicable – no new data generated.

\small{

$\noindent$

School of Mathematics and Statistics

University of Sheffield

Hicks Building

Sheffield S3 7RH

UK

email: v.bavula@sheffield.ac.uk

\end{document}